\documentclass[12pt,a4paper]{amsart}
\usepackage{t1enc}
\usepackage[latin1]{inputenc}
\usepackage[english]{babel}
\usepackage{hyperref}
\usepackage{graphicx}

\usepackage{latexsym}
\usepackage{amsmath,amsfonts,amssymb,amsthm}
\usepackage{mathptmx}

\begin{document}

\title {
	Contact geometry of one dimensional holomorphic foliations
}

\author{Giuseppe Tomassini}
\address{
	G. Tomassini:
	Scuola Normale Superiore,
	Piazza dei Cavalieri,
	7 --- I-56126 Pisa, Italy}
\email{g.tomassini@sns.it}

\author{Sergio Venturini}
\address{
	S. Venturini:
	Dipartimento Di Matematica,
	Universit\`{a} di Bologna,
	\,\,Piazza di Porta S. Donato 5 ---I-40127 Bologna,
	Italy}
\email{venturin@dm.unibo.it}


%
%
\def\R{{\rm I\kern-.185em R}}
\def\C{{\rm\kern.37em\vrule height1.4ex width.05em depth-.011em\kern-.37em C}}
\def\N{{\rm I\kern-.185em N}}
\def\Z{{\bf Z}}
\def\Q{{\mathchoice{\hbox{\rm\kern.37em\vrule height1.4ex width.05em 
depth-.011em\kern-.37em Q}}{\hbox{\rm\kern.37em\vrule height1.4ex width.05em 
depth-.011em\kern-.37em Q}}{\hbox{\sevenrm\kern.37em\vrule height1.3ex 
width.05em depth-.02em\kern-.3em Q}}{\hbox{\sevenrm\kern.37em\vrule height1.3ex
width.05em depth-.02em\kern-.3em Q}}}}
\def\P{{\rm I\kern-.185em P}}
\def\H{{\rm I\kern-.185em H}}
%
\def\Aleph{\aleph_0}
\def\ALEPH#1{\aleph_{#1}}
\def\sset{\subset}\def\ssset{\sset\sset}
%
\def\bar#1{\overline{#1}}
\def\dim{\mathop{\rm dim}\nolimits}
\def\half{\textstyle{1\over2}}
\def\Half{\displaystyle{1\over2}}
\def\mlog{\mathop{\half\log}\nolimits}
\def\Mlog{\mathop{\Half\log}\nolimits}
\def\Det{\mathop{\rm Det}\nolimits}
\def\Hol{\mathop{\rm Hol}\nolimits}
\def\Aut{\mathop{\rm Aut}\nolimits}
\def\Re{\mathop{\rm Re}\nolimits}
\def\Im{\mathop{\rm Im}\nolimits}
\def\Ker{\mathop{\rm Ker}\nolimits}
\def\Fix{\mathop{\rm Fix}\nolimits}
\def\Res{\mathop{\rm Res}\nolimits}
\def\sp{\mathop{\rm sp}\nolimits}
\def\id{\mathop{\rm id}\nolimits}
\def\Trace{\mathop{\rm Tr}\nolimits}
\def\cancel#1#2{\ooalign{$\hfil#1/\hfil$\crcr$#1#2$}}
\def\prevoid{\mathrel{\scriptstyle\bigcirc}}
\def\void{\mathord{\mathpalette\cancel{\mathrel{\scriptstyle\bigcirc}}}}
\def\n{{}|{}\!{}|{}\!{}|{}}
\def\abs#1{\left|#1\right|}
\def\norm#1{\left|\!\left|#1\right|\!\right|}
\def\nnorm#1{\left|\!\left|\!\left|#1\right|\!\right|\!\right|}
%
\def\upperint{\int^{{\displaystyle{}^*}}}
\def\lowerint{\int_{{\displaystyle{}_*}}}
\def\Upperint#1#2{\int_{#1}^{{\displaystyle{}^*}#2}}
\def\Lowerint#1#2{\int_{{\displaystyle{}_*}#1}^{#2}}
%
\def\rem #1::#2\par{\medbreak\noindent{\bf #1}\ #2\medbreak}
\def\proclaim #1::#2\par{\removelastskip\medskip\goodbreak{\bf#1:}
\ {\sl#2}\medskip\goodbreak}
\def\ass#1{{\rm(\rmnum#1)}}
\def\assertion #1:{\Acapo\llap{$(\rmnum#1)$}$\,$}
\def\Assertion #1:{\Acapo\llap{(#1)$\,$}}
\def\acapo{\hfill\break\noindent}
\def\Acapo{\hfill\break\indent}
\def\proof{\removelastskip\par\medskip\goodbreak\noindent{\it Proof.\/\ }}
\def\prova{\removelastskip\par\medskip\goodbreak
\noindent{\it Dimostrazione.\/\ }}
\def\qed{{\bf //}\par\smallskip}
\def\BeginItalic#1{\removelastskip\par\medskip\goodbreak
\noindent{\it #1.\/\ }}
\def\iff{if, and only if,\ }
\def\sse{se, e solo se,\ }
\def\rmnum#1{\romannumeral#1{}}
\def\Rmnum#1{\uppercase\expandafter{\romannumeral#1}{}}
\def\smallfrac#1/#2{\leavevmode\kern.1em
\raise.5ex\hbox{\the\scriptfont0 #1}\kern-.1em
/\kern-.15em\lower.25ex\hbox{\the\scriptfont0 #2}}
%
\def\Left#1{\left#1\left.}
\def\Right#1{\right.^{\llap{\sevenrm
\phantom{*}}}_{\llap{\sevenrm\phantom{*}}}\right#1}
%
%
%
\def\dimens{3em}
\def\symb[#1]{\noindent\rlap{[#1]}\hbox to \dimens{}\hangindent=\dimens}
\def\references{\bigskip\noindent{\bf References.}\bigskip}
\def\art #1 : #2 ; #3 ; #4 ; #5 ; #6. \par{#1, 
{\sl#2}, #3, {\bf#4}, (#5), #6.\par\smallskip}
\def\book #1 : #2 ; #3 ; #4. \par{#1, {\bf#2}, #3, #4.\par\smallskip}
\def\freeart #1 : #2 ; #3. \par{#1, {\sl#2}, #3.\par\smallskip}
%
%
%
%
\def\name{\hbox{Sergio Venturini}}
\def\snsaddress{\indent
\vbox{\bigskip\bigskip\bigskip
\name
\hbox{Scuola Normale Superiore}
\hbox{Piazza dei Cavalieri, 7}
\hbox{56126 Pisa (ITALY)}
\hbox{FAX 050/563513}}}
\def\cassinoaddress{\indent
\vbox{\bigskip\bigskip\bigskip
\name
\hbox{Universit\`a di Cassino}
\hbox{via Zamosch 43}
\hbox{03043 Cassino (FR)}
\hbox{ITALY}}}
\def\bolognaaddress{\indent
\vbox{\bigskip\bigskip\bigskip
\name
\hbox{Dipartimento di Matematica}
\hbox{Universit\`a di Bologna}
\hbox{Piazza di Porta S. Donato 5}
\hbox{40127 Bologna (BO)}
\hbox{ITALY}
\hbox{venturin@dm.unibo.it}
}}
\def\homeaddress{\indent
\vbox{\bigskip\bigskip\bigskip
\name
\hbox{via Garibaldi, 7}
\hbox{56124 Pisa (ITALY)}}}
\def\doubleaddress{
\vbox{
\hbox{\name}
\hbox{Universit\`a di Cassino}
\hbox{via Zamosch 43}
\hbox{03043 Cassino (FR)}
\hbox{ITALY}
\smallskip
\hbox{and}
\smallskip
\hbox{Scuola Normale Superiore}
\hbox{Piazza dei Cavalieri, 7}
\hbox{56126 Pisa (ITALY)}
\hbox{FAX 050/563513}}}
\def\sergio{{\rm\bigskip
\centerline{Sergio Venturini}
\leftline{\bolognaaddress}
\bigskip}}
%
%
%
%
%

\newtheorem{theorem}{Theorem}[section]
\newtheorem{proposition}{Proposition}[section]
\newtheorem{lemma}{Lemma}[section]
\newtheorem{corollary}{Corollary}[section]
\newtheorem{remark}{Remark}[section]
\newtheorem{definition}{Definition}[section]

\newtheorem{teorema}{Teorema}[section]
\newtheorem{proposizione}{Proposizione}[section]
\newtheorem{corollario}{Corollario}[section]
\newtheorem{osservazione}{Osservazione}[section]
\newtheorem{definizione}{Definizione}[section]
\newtheorem{esempio}{Esempio}[section]
\newtheorem{esercizio}{Esercizio}[section]
\newtheorem{congettura}{Congettura}[section]

\bibliographystyle{abbrv}

\def\CMan{M}
\def\dimCMan{n}
\def\dd{{\rm d}}
\def\dc{{\dd^c}}
\def\VField{\xi}
\def\Eq{u}
\def\Leaf{S}
\def\CLocus{Z}
\def\CForm{\theta}
\def\HyperSurface{V}
\def\Point{p}
\def\ia{i}
\def\ib{j}
\def\ic{m}
\def\Smooth{k}
\def\Inclusion{j}
\def\NeighbA{{\CMan_0}}
\def\NeighbB{{\CMan_1}}
\def\NeighbC{{\CMan_2}}
\def\ConstA{c}
\def\res{\mathop{\hbox{\vrule height 7pt width .5pt depth 0pt\vrule height .5pt width 6pt depth 0pt}}\nolimits}

\begin{abstract}
Let $V$ be a real hypersurface of class ${\rm C}^k$, $k\ge 3$, in a complex manifold $M$ of complex dimension $n+1$, $HT(V)$ the holomorphic tangent bundle to $V$ giving the induced {\rm CR} structure on $V$. Let $\theta$ be a contact form for $(V,HT(V))$, $\xi_0$ the Reeb vector field determined by $\theta$ and assume that  $\xi_0$ is of class ${\rm C}^k$. In this paper we prove the following theorem (cf. Theorem \ref{thm::CauchyMA}): if the integral curves of $\xi_0$ are real analytic then there exist an open neighbourhood $\NeighbA\sset\CMan$ of $\HyperSurface$ and a solution $\Eq\in C^\Smooth(\NeighbA)$ of the complex Monge-Ampère equation $(\dd\dc\Eq)^{\dimCMan+1}=0$ on $\NeighbA$ which is a defining equation for $V$. Moreover, the Monge-Ampère foliation associated to $u$ induces on $V$ that one associated to the Reeb vector field.
The converse is also true.
The result is obtained solving a Cauchy problem for infinitesimal symmetries
of CR distributions of codimension one which is of independent interest
(cf. Theorem \ref{thm::CauchyCF} below).

\end{abstract}

\maketitle

\nocite{book:KobayashiNomizuA}
\nocite{book:KobayashiNomizuB}
\nocite{book:KushnerLychaginRubtsov}
\nocite{book:GilbertBuchanan}
\nocite{book:BlairContactEtc}
\nocite{article:AndreottiFredricksEmbeddingCR}
\nocite{article:BedfordKalka}
\nocite{article:DuchampKalka}
\nocite{article:BedfordBurns}


\section{\label{section:SectionIntro}Introduction}

We follow \cite{book:KobayashiNomizuA}, \cite{book:KobayashiNomizuB}
for standard notations in differential geometry and
complex manifolds.

Let $\CMan$ be a complex manifold of complex dimension $\dimCMan+1$
with complex structure $J$ on the tangent bundle  $T(\CMan)$.

A function $\Eq\in C^2(\CMan)$ is a solution of the
\emph{complex Monge-Amp\`ere equation} if and only if 
$$
MA(\Eq)=(\dd\dc\Eq)^{\dimCMan+1}
	=(2i\partial\overline\partial\Eq)^{\dimCMan+1}=0
$$
where $\dc=i(\overline\partial-\partial)$;
in local holomorphic coordinates $z_1,\dots,z_{n+1}$  
$$
\det\left(\frac{\partial^2 \Eq}{\partial z_\alpha\overline z_\beta}\right)=0.
$$

If $\omega$ is a $1-$differential form,
then we denote by $\omega^c$ the $1-$differential form
which satisfies $\omega^c(X)=-\omega(JX)$,
for each vector field $X$,
so that if $f\in C^1(\CMan)$ then $(\dd f)^c=\dc f$.


If $\HyperSurface\subset\CMan$ is a real hypersurface of class $C^\Smooth$, $k>1$,
set 
$$
HT(\HyperSurface)=\{X\in T(\HyperSurface)\mid JX\in T(\HyperSurface)\}.
$$
Then $HT(\HyperSurface)$ is a $J$-invariant distribution 
of real dimension $2\dimCMan$, called the \emph{holomorphic tangent bundle} to $\HyperSurface$,
The distribution $HT(\HyperSurface)$ gives $V$ a {\rm CR} structure, that induced by the complex structure of $M$; $V$ endowed with this {\rm CR} structure is denoted by $(\HyperSurface, HT(\HyperSurface))$

We recall that the Levi form of $\HyperSurface$ is non degenerate if, and only if, locally there exists
a real differential form $\CForm$ of degree $1$ on $\HyperSurface$ such that
the restriction to $HT(\HyperSurface)$ of the skew 2-form $\dd\CForm$ is non degenerate.
In this case $(\HyperSurface, HT(\HyperSurface))$ is a contact manifold with a contact form  $\CForm$, also said a \emph{ contact ${\rm CR}$ hypersurface}.

We refer to \cite{book:BlairContactEtc} and \cite{book:KushnerLychaginRubtsov}
for basic facts on contact geometry.

A vector field $\VField$ on the contact manifold $(\HyperSurface, HT(\HyperSurface))$
is an \emph{infinitesimal symmetry} if $[\VField,HT(\HyperSurface)]\sset HT(\HyperSurface)$;
an infinitesimal symmetry such that $\VField(\Point)\in HT_\Point(\HyperSurface)$
for each point $\Point\in\HyperSurface$
is a \emph{characteristic vector field} for the distribution HT(\HyperSurface)
(see e. g. Section 1.2 of \cite{book:KushnerLychaginRubtsov}).

If $\CForm$ is a contact form for the
contact manifold $(\HyperSurface, HT(\HyperSurface))$ then there exists
a unique vector field $\VField_0$ on $\HyperSurface$,
the \emph{Reeb vector field}, which satisfies
$\CForm(\VField_0)=1$ and $\dd\CForm(\VField_0,X)=0$ for each
vector field $X$ on $\HyperSurface$.
It is easy to show that $\VField_0$ is an infinitesimal symmetry
of the distribution $HT(\HyperSurface)$.

Let $\HyperSurface\sset\CMan$ be a hypersurface of class $C^\Smooth$.
A real function $\Eq\in C^k(\CMan)$ is called an \emph{equation}
of $\HyperSurface\sset\CMan$ if $\HyperSurface=\{\Eq=0\}$ and $\dd\Eq\neq 0$ near $\HyperSurface$.
Then the restriction of the form $\dc\Eq$ to $T(\HyperSurface)$ is a real $1$-form
which defines the {\rm CR} structure $HT(\HyperSurface)$.
Observe that $\HyperSurface$ is a contact ${\rm CR} $ hypersurface if, and only if,
the $2(\dimCMan+1)$ form $\dd\Eq\wedge \dc\Eq\wedge(\dd\dc\Eq)^{\dimCMan}$ does not
vanish on (a neighbourhood of) $\HyperSurface$.

In this paper we are interested in studying the existence of equations
$\Eq\in C^\Smooth(\CMan)$ of a given hypersurface
$\HyperSurface\sset\CMan$ wich are solution of the
complex Monge-Amp\`ere equation $(\dd\dc\Eq)^{\dimCMan+1}=0$ (where $\dimCMan+1$ is the complex dimension of $\CMan$). 
We prove the following theorem (cf. Theorem \ref{thm::CauchyMA}): if the integral curves of $\xi_0$ are real analytic then there exist an open neighbourhood $\NeighbA\sset\CMan$ of $\HyperSurface$ and a solution $\Eq\in C^\Smooth(\NeighbA)$ of the complex Monge-Amp\`ere equation $(\dd\dc\Eq)^{\dimCMan+1}=0$ on $\NeighbA$ which is a defining equation for $V$. Moreover, the Monge-Amp\`ere foliation associated to $u$ induces on $V$ that one associated to the Reeb vector field.
The converse is also true.

The result is obtained solving a Cauchy problem for infinitesimal symmetries
(cf. Theorem \ref{thm::CauchyCF} below)
of CR distributions of codimension one which is of independent interest

As for the contents of the paper, in Section \ref{section:VectorFields} we define the notion of
\emph{calibrated foliation} which  is nothing but than a pair $(\VField,\Eq)$
where $\VField$ is a vector field on the complex manifold $\CMan$ such that
$[\VField,J\VField]=0$ and $\Eq$ is a function on $\CMan$ which satisfy
$\dc\Eq(\VField)=0$ and $\dd\Eq(\VField)=1$.
If $(\VField,\Eq)$ is calibrated foliation then the vector field $\VField$
induces on $\CMan$ a holomorphic foliation whose leaves are Riemann surfaces.
The main result of the section is, roughly speaking, that the set $\CLocus$
of the points of the complex manifold $\CMan$ where the vector field $\VField$
is an infinitesimal simmetry for the distribution ${\sf Ker}\,\dd\Eq\cap{\sf Ker}\,\dc\Eq$
intersects each leaf $\Leaf$ of the holomorphic foliation along an analytic subset
of $\Leaf$; hence either $\Leaf\sset\CLocus$ or $\Leaf\cap\CLocus$ is a discrete
subset of $\Leaf$ (cf. Theorem \ref{thm:main}).
Here the basic tools are provided by the theory of the
\emph{generalized analytic functions}, developped extensively in \cite{book:GilbertBuchanan}, dealing with functions wich satisfy a 
first-order complex linear differential system of equations having the Cauchy-Riemann
operator as principal part symbol (cf. Theorem \ref{thm::ZeroesHyperanalitic}).

In Section \ref{section:SectionCauchy} we prove that if $\HyperSurface\sset\CMan$
is a real hypersurface of the complex manifold $\CMan$ and $\VField_0$ is a
vector field on $\HyperSurface$ having real analytic integral curves 
then there exists locally a unique calibrated foliation $(\VField,\Eq)$
in a neighbourhood of $\HyperSurface$ in $\CMan$
such that $\Eq$ vanishes on $\HyperSurface$ and $\VField$ extends $\VField_0$.
The results of the previous section are used in order to show that the vector field $\VField$
is an infinitesimal simmetry of the distribution ${\sf Ker}\,\dd\Eq\cap{\sf Ker}\,\dc\Eq$
if, and only, if $\VField_0$ is an infinitesimal simmetry of the distribution
$HT(\HyperSurface)$ (cf. Theorem \ref{thm::CauchyCF}).

If $(\VField,\Eq)$ is a calibrated foliation and $\VField$ is an
infinitesimal simmetry of the distribution ${\sf Ker}\,\dd\Eq\cap{\sf Ker}\,\dc\Eq$
then $\Eq$ is a solution of the complex Monge-Amp\`ere equation on $\CMan$.
This elementary observation yields to the main result of Section 
\ref{section:ComplexMA} on the existence, for a Levi non degenerate CR real hypersurface $\HyperSurface$ with assigned contact form $\CForm$, of an equation which is solution of the complex Monge-Amp\`ere. The result is that such equation exists if, and only if, the Reeb vector field $\VField_0$ associated to the contact manifold $(\HyperSurface,\CForm)$ has real analytic integral curves (cf. Theorem \ref{thm::CauchyMA}). Indeed in this case by the results of the
previous section there exists a unique calibrated foliation $(\VField,\Eq)$
in a neighbourhood of $\HyperSurface$ such that $\VField$ extends the
Reeb vector field $\VField_0$; then, since the Reeb vector field is an infinitesimal simmetry of the corresponding contact structure, $\Eq$ is a solution of the complex Monge-Amp\`ere equation.

Finally, let us observe that, by a result of Andreotti and Fredricks
(cf. \cite{article:AndreottiFredricksEmbeddingCR}, Theorem 1.12)
any real analytic codimension one CR manifold with
integrable real analytic CR distribution embeds
in some complex manifold as a CR hypersurface. Thus, our theory applies to these abstract CR manifolds too.

\section{\label{section:VectorFields}Calibrated one dimensional foliations}

Let $\CMan$ be a complex manifold of dimension $\dimCMan+1$
 with (integrable) complex structure $J$.

Let $\VField$ be a vector field on $\CMan$
and let $\Eq$ be a funcion on $\CMan$. 
We say that the pair $(\VField,\Eq)$ is a
\emph{calibrated foliation of dimension one}
(or simply a calibrated foliation)
of class $C^\Smooth$, $\Smooth>0$, if
$\VField$ and $\Eq$ are of class $C^\Smooth$ and satisfies
the conditions
\begin{eqnarray}
\label{eq:CalibFoliationA}
&&[\VField,J\VField]=0,
\\
\label{eq:CalibFoliationB}
&&\dd\Eq(\VField)=0,\quad \dc\Eq(\VField)=1.
\end{eqnarray}

From $\dc\Eq(\VField)=1$ it follows that
$\dc\Eq$ never vanishes identically at any point $\Point\in\CMan$
and hence for each constant $c\in\R$ the subset of points
$\Point\in\CMan$ which satisfie
$\Eq(\Point)=\ConstA$ is either the empty set or is a real
hypersurface of $\CMan$ of class $C^\Smooth$.

Clearly the vector fields $\VField$ and $J\VField$
generate a $J$-invariant bidimensional distribution on $TM$
whoose maximal (connected) integral submanifolds are Riemann surfaces
which fill the manifol $\CMan$.
We call such Riemann surfaces the \emph{leafs} of the 
calibrated foliation $(\VField,\Eq)$.

Let $\Leaf$ be a leaf of the calibrated foliation $(\VField,\Eq)$.
An \emph{adapded holomorphic coordinate} on $\Leaf$ is
a holomorphic map $z:A\to\C$ where $A\sset\Leaf$ is open in $\Leaf$
and $\Im z=\Eq_{|A}$.
It is straighforward to prove that for each $\Point\in\Leaf$ there exists
an adapded holomorphic coordinate on $\Leaf$ defined in a neighbourhood
of $\Point$ in $\Leaf$.
This shows in particular that the restriction of the function $\Eq$
to each leaf $\Leaf$ is a harmonic function on $\Leaf$.

If $(\VField,\Eq)$ is a calibrated foliation then ${\sf Ker}\, \dd\Eq\sset T(\CMan)$ is
an integrable distribution having as maximal integrable submanifold
the connected components of the hypersurfaces defined by as $\Eq=\ConstA$,
$c$ real constant.
Let us observe that if $V$ is a
(maximal) integral submanifold of ${\sf Ker}\,  \dd\Eq$
and if $\Leaf$ is a leaf of $(\VField,\Eq)$
then $\HyperSurface\cap\Leaf$ is a
(maximal) integral curve of the vector field $\VField$,
and each (maximal) integral curve of the vector field $\VField$
is obtained in this way.

Let $(\VField,\Eq)$ be a calibrated foliation.
We define the \emph{contact locus} of the calibrated foliation $(\VField,\Eq)$
as the set $\CLocus$ of the points $p\in\CMan$ where the
differential form $L_\VField(\dc\Eq)$ vanishes.
Here $L_\VField$ stands for the Lie derivative with respect
to the vector field $\VField$.
We also say that $(\VField,\Eq)$ is a \emph{contact calibrated foliation}
if $\CLocus=\CMan$.
It is easy to show that $\VField$ is a characteristic vector field of the
distribution ${\sf Ker}\,  d\Eq$ and that
$(\VField,\Eq)$ is a contact calibrated foliation if, and only if
the vector field $\VField$ is an infinitesimal symmetry for the
distribution ${\sf Ker}\,\dd\Eq\cap{\sf Ker}\, \dc\Eq$
(cf. e. g. Theorem 1.2.1 of \cite{book:KushnerLychaginRubtsov}).

The main result of this section is the following:

\begin{theorem}\label{thm:main}
Let $(\VField,\Eq)$ be a calibrated foliation of class $C^\Smooth$,
$\Smooth\geq2$, on the complex manifold $\CMan$ with contact locus $\CLocus$.

If $\Leaf$ is a leaf of the calibrated foliation $(\VField,\Eq)$
then $\Leaf\cap\CLocus$ is an analytic subset of $\Leaf$,
that is $\Leaf\cap\CLocus\sset\Leaf$ is locally defined by the common
zeroes of holomorphic function on $\Leaf$ and hence
either $\Leaf\sset\CLocus$ or $\Leaf\cap\CLocus$ is a
discrete subset of $\Leaf$.
\end{theorem}

Before getting involved in the  proof let us observe an immediate consequence of Theorem 3.12 of \cite{book:GilbertBuchanan}:
 
\begin{theorem}\label{thm::ZeroesHyperanalitic}
Let $D\sset\C$ be an open domain and let
$w_\ia\in C^1(D),\ia=1,\ldots,p$
and $A_{\ia\ib},B_{\ia\ib}\in C^0(D),\ia,\ib=1,\ldots,p$
be complex functions. Assume that
\begin{equation}\label{eq::PascaliSystem}
	\frac{\partial w_\ia}{\partial\bar z}=
	\sum_{\ib=1}^p A_{\ia\ib}w_\ib+B_{\ia\ib}\bar w_\ib
	\quad \ia=1,\ldots,p
	\end{equation}
holds on $D$.
Then the common zeroes of the functions $w_\ia$ is an analytic
subset of $D$ and hence
either the functions $w_\ia$ vanishes identically on $D$
or the common zeroes of the functions $w_\ia$ is a discrete set of $D$.
\end{theorem}

Solutions of (\ref{eq::PascaliSystem}) are called \emph{generalized analytic functions}.
We also summarize by the following lemma some elementary fact that will be used in the sequel.
\begin{lemma}\label{lemma::ElementaryLemma}
Let $\CMan$ a complex manifold
and let $X,Y$ be vector fields on $\CMan$

If $f\in C^2(\CMan)$ then
	\begin{eqnarray}\label{eq::lemma::Jinvariance}
		&&\dd\dc f(JX, JY)=\dd\dc f(X, Y).
	\end{eqnarray}

If $\theta$ is a differential form on $M$ of degree one
	and $\theta(X)=c_1$, $\theta(Y)=c_2$ with $c_1$ and $c_2$ constant
	then 
	\begin{eqnarray}\label{eq::lemma::b}
		&&\dd\theta(X, Y)=(L_X\theta)(Y)=-\theta\bigl([X, Y]\bigr).
	\end{eqnarray}
	
If $f\in C^2(\CMan)$
	and $\dd f(X)=c_1$, $\dd f(Y)=c_2$ with $c_1$ and $c_2$ constant
	then 
	\begin{eqnarray}\label{eq::lemma::c}
		&&\dd f\bigl([X, Y]\bigr)=0.
	\end{eqnarray}
\end{lemma}

{\it Proof of theorem \ref{thm:main}.\/\ }
We will prove the theorem showing that for each
leaf $\Leaf$ the set $\Leaf\cap\CLocus$ is locally
the set of common zeroes of 
a set of funcion $w_1,\ldots,w_\dimCMan$ which satisfy a system
of equations of the form (\ref{eq::PascaliSystem}).
Let $\dimCMan+1$ be the dimension of $\CMan$.
Let $\Leaf$ be a leaf.
Let $\Point\in\Leaf$.
For a neighbourhood $U$ of $\Point$ in $\CMan$ small enought
there exist an adapted holomorphic coordinate $z:\Leaf\cap U\to\C$
and $C^2$ vector fields $X_1,\ldots,X_\dimCMan$ on $U$ such that
$\VField,X_1,JX_1,\ldots,X_\dimCMan,JX_\dimCMan$ generate the
distribution ${\sf Ker}\,  \dd\Eq$.

Of course $X_1,JX_1,\ldots,X_\dimCMan,JX_\dimCMan$ generate the
distribution ${\sf Ker}\,  \dd\Eq\cap {\sf Ker}\,  \dc\Eq$ and 
$\VField,J\VField,X_1,JX_1,\ldots,X_\dimCMan,JX_\dimCMan$ generate the
whole $T(\CMan)$.

Setting $\omega=L_\VField \dc(\Eq)$ and for $\ia=1,\ldots,\dimCMan$
let define $u_\ia=\omega(X_\ia)$, $v_\ia=\omega^c(X_\ia)$.

Since $\omega(\VField)=\omega(J\VField)=0$ it follows that
$Z\cap U$ is exactly the common zero set of the functions
$u_\ia,v_\ia$, $\ia=1,\ldots,\dimCMan$.

We now first prove that for $\ia=1,\ldots,\dimCMan$,
\begin{eqnarray}
\label{eq:main:DefA}
&&\dc\Eq\bigl([X_\ia,\VField]\bigr)=u_\ia,
\\
\label{eq:main:DefB}
&&\dc\Eq\bigl([JX_\ia,\VField]\bigr)=v_\ia,
\\
\label{eq:main:DefC}
&&\dc\Eq\bigl([X_\ia,J\VField]\bigr)=-v_\ia,
\\
\label{eq:main:DefD}
&&\dc\Eq\bigl([JX_\ia,J\VField]\bigr)=u_\ia.
\end{eqnarray}

Indeed, since $\dc\Eq(X_\ia)=\dc\Eq(JX_\ia)=0$,
\begin{eqnarray}
\label{eq::nolabel::a}\nonumber
&&u_\ia
	=L_\VField \dc\Eq(X_\ia)
	=\VField\bigl(\dc\Eq(X_\ia)\bigr)
		-\dc\Eq\bigl([\VField,X_\ia]\bigr)
	=\dc\Eq\bigl([X_\ia,\VField]\bigr).
\end{eqnarray}
and
\begin{eqnarray}
\label{eq::nolabel::b}\nonumber
&&v_\ia
	=L_\VField \dc\Eq(JX_\ia)
	=\VField\bigl(\dc\Eq(JX_\ia)\bigr)
		-\dc\Eq\bigl([\VField,JX_\ia]\bigr)
	=\dc\Eq\bigl([JX_\ia,\VField]\bigr),
\end{eqnarray}
which proves (\ref{eq:main:DefA}) and (\ref{eq:main:DefB}).

Using (\ref{eq::lemma::Jinvariance}), (\ref{eq::lemma::b})
and the identity $J^2(X)=-X$ we obtain
\begin{eqnarray}
\nonumber
	\dc\Eq\bigl([X_\ia,J\VField]\bigr)
	&=&-\dd\dc\Eq(X_\ia,J\VField)
	=\dd\dc\Eq(JX_\ia,\VField)
\\
\nonumber
	&=&-\dc\Eq\bigl([JX_\ia,\VField]\bigr)
	=-v_\ia.
\end{eqnarray}
and
\begin{eqnarray}
\nonumber
	\dc\Eq\bigl([JX_\ia,J\VField]\bigr)
	&=&-\dd\dc\Eq(JX_\ia,J\VField)
	=-\dd\dc\Eq(X_\ia,\VField)
\\
\nonumber
	&=&\dc\Eq\bigl([X_\ia,\VField]\bigr)
	=u_\ia.
\end{eqnarray}
which proves (\ref{eq:main:DefC}) and (\ref{eq:main:DefD}).


Then we prove that for $\ia=1,\ldots,\dimCMan$
\begin{eqnarray}
\label{eq::commutator::a}
	[X_\ia,\VField]&=&u_\ia\VField + 
		\sum_{\ib=1}^\dimCMan\bigl(
			a_{\ia\ib}X_\ib+b_{\ia\ib}JX_\ib
		\bigr),
\\
\label{eq::commutator::b}
	[JX_\ia,\VField]&=&v_\ia\VField + 
		\sum_{\ib=1}^\dimCMan\bigl(
			c_{\ia\ib}X_\ib+d_{\ia\ib}JX_\ib
		\bigr),
\\
\label{eq::commutator::c}
	[X_\ia,J\VField]&=&-v_\ia\VField + 
		\sum_{\ib=1}^\dimCMan\bigl(
			e_{\ia\ib}X_\ib+f_{\ia\ib}JX_\ib
		\bigr),
\\
\label{eq::commutator::d}
	[JX_\ia,J\VField]&=&u_\ia\VField + 
		\sum_{\ib=1}^\dimCMan\bigl(
			g_{\ia\ib}X_\ib+h_{\ia\ib}JX_\ib
		\bigr),
\end{eqnarray}
where $a_{\ia\ib},\ldots,h_{\ia\ib}$ are $C^1$ functions on $U$.

Indeed consider first $[X_\ia,\VField]$.
Since $\dc\Eq(X_\ia)=0, \dc\Eq(\VField)=0$ by (\ref{eq::lemma::c})
it follows that $[X_\ia,\VField]\in{\sf Ker}\, \dd\Eq$ and hence

\begin{eqnarray}
\label{eq::commutator::nonumbera}
	[X_\ia,\VField]&=&\lambda_\ia\VField + 
		\sum_{\ib=1}^\dimCMan\bigl(
			a_{\ia\ib}X_\ib+b_{\ia\ib}JX_\ib
		\bigr)
\end{eqnarray}
where $\lambda_\ia$, $a_{\ia\ib}$ and $b_{\ia\ib}$
are some $C^1$ functions on $U$.

Since $\dc\Eq(\VField)=1$ and
$\dc\Eq(X_\ia)=\dc\Eq(JX_\ia)=0$ it follows that
\begin{eqnarray}
\label{eq::commutator::nonumberb}
	&&u_\ia=\dc\Eq\bigl([X_\ia,\VField]\bigr)=\lambda_\ia,
\end{eqnarray}
and this proves (\ref{eq::commutator::a}).

The proofs of (\ref{eq::commutator::b}),
(\ref{eq::commutator::c}) and (\ref{eq::commutator::d})
are similar.

Finally we claim that
\begin{eqnarray}
\label{eq::bicommutator::a}
	\dc\Eq\bigl([[X_\ia,\VField],J\VField]\bigr)
		&=&-J\VField(u_\ia)+\sum_{\ib=1}^\dimCMan
			\bigl(b_{\ia\ib}u_\ib-a_{\ia\ib}v_\ib\bigr),
	\\
\label{eq::bicommutator::b}
	\dc\Eq\bigl([[JX_\ia,\VField],J\VField]\bigr)
		&=&-J\VField(v_\ia)+\sum_{\ib=1}^\dimCMan
			\bigl(d_{\ia\ib}u_\ib-c_{\ia\ib}v_\ib\bigr),
	\\
\label{eq::bicommutator::c}
	\dc\Eq\bigl([[X_\ia,J\VField],\VField]\bigr)
		&=&\phantom{-J}\VField(v_\ia)+\sum_{\ib=1}^\dimCMan
			\bigl(e_{\ia\ib}u_\ib+f_{\ia\ib}v_\ib\bigr).
	\\
\label{eq::bicommutator::d}
	\dc\Eq\bigl([[JX_\ia,J\VField],\VField]\bigr)
		&=&\phantom{J}-\VField(u_\ia)+\sum_{\ib=1}^\dimCMan
			\bigl(g_{\ia\ib}u_\ib+h_{\ia\ib}v_\ib\bigr).
\end{eqnarray}

Indeed from (\ref{eq::commutator::a}) we obtain
\begin{eqnarray}
	\label{eq::bicommutator::nonumbera}
	[[X_\ia,\VField],J\VField]
	&=&[u_\ia\VField + 
		\sum_{\ib=1}^\dimCMan\bigl(
			a_{\ia\ib}X_\ib+b_{\ia\ib}JX_\ib
		\bigr),J\VField]
	\nonumber\\
	\label{eq::bicommutator::nonumberb}
	&=&-J\VField(u_\ia)\VField +
		\sum_{\ib=1}^\dimCMan\bigl(
			a_{\ia\ib}[X_\ib,J\VField]+b_{\ia\ib}[JX_\ib,J\VField]
		\bigr)
	\nonumber\\
		&&-\sum_{\ib=1}^\dimCMan\bigl(
			J\VField(a_{\ia\ib})X_\ib+J\VField(b_{\ia\ib})JX_\ib
		\bigr).
	\nonumber
\end{eqnarray}
Applying $\dc\Eq$, using (\ref{eq:main:DefC}) and (\ref{eq:main:DefD})
we obtain
\begin{eqnarray}
	\label{eq::bicommutator::nonumberc}
		\dc\Eq\bigl([[X_\ia,\VField],J\VField]\bigr)
		&=&-J\VField(u_\ia)
	\nonumber\\
	\label{eq::bicommutator::nonumberd}
			&&+\sum_{\ib=1}^\dimCMan\Bigl(
				a_{\ia\ib}\dc\Eq\bigl([X_\ib,J\VField]\bigr)
				+b_{\ia\ib}\dc\Eq\bigl([JX_\ib,J\VField]\bigr)
			\Bigr)
	\nonumber\\
	\label{eq::bicommutator::nonumbere}
		&=&-J\VField(u_\ia)
			+\sum_{\ib=1}^\dimCMan\bigl(
				b_{\ia\ib}u_\ib
				-a_{\ia\ib}v_\ib
			\bigr),
	\nonumber
\end{eqnarray}
and this proves (\ref{eq::bicommutator::a}).

The proofs of (\ref{eq::bicommutator::b}),
(\ref{eq::bicommutator::c}) and (\ref{eq::bicommutator::d})
are similar.

The relation $[\VField,J\VField]=0$ and the Jacobi identity for the
Poisson bracket yield
\begin{eqnarray}
	\label{eq::poisson::a}\nonumber
	[[JX_\ia,\VField],J\VField]&=&[[JX_\ia,J\VField],\VField],
	\\
	\label{eq::poisson::b}\nonumber
	[[X_\ia,\VField],J\VField]&=&[[X_\ia,J\VField],\VField].
\end{eqnarray}
Applying $\dc\Eq$ and using
(\ref{eq::bicommutator::a}),$\ldots$,(\ref{eq::bicommutator::d}),
after some rearrangement we obtain that for $\ia=1,\ldots,\dimCMan$
\begin{eqnarray}
	\label{eq::precr::a}
	\VField(u_\ia)-J\VField(v_\ia)&=&
		\sum_{\ib=1}^\dimCMan\bigl[
			(g_{\ia\ib}-d_{\ia\ib})u_\ib+
			(h_{\ia\ib}+c_{\ia\ib})v_\ib
		\bigr],
	\\
	\label{eq::precr::b}
	J\VField(u_\ia)+\VField(v_\ia)&=&
		\sum_{\ib=1}^\dimCMan\bigl[
			(b_{\ia\ib}-e_{\ia\ib})u_\ib-
			(a_{\ia\ib}+f_{\ia\ib})v_\ib
		\bigr].
\end{eqnarray}

Considering the holomorphic coordinate $z:\Leaf\cap U$,
if $z=x+y$ with $x$ and $y$ real functions we have
\begin{eqnarray}
	\label{eq::cr::nonumbera}\nonumber
	\VField&=&\frac{\partial}{\partial x},
	\\
	\label{eq::cr::nonumberb}\nonumber
	J\VField&=&\frac{\partial}{\partial y},
\end{eqnarray}
and hence in such coordinate system the equations 
(\ref{eq::precr::a}) and (\ref{eq::precr::b})
reduce to
\begin{eqnarray}
	\label{eq::cr::a}
	\frac{\partial u_\ia}{\partial x}-\frac{\partial v_\ia}{\partial y}&=&
		\sum_{\ib=1}^\dimCMan\bigl[
			(g_{\ia\ib}-d_{\ia\ib})u_\ib+
			(h_{\ia\ib}+c_{\ia\ib})v_\ib
		\bigr],
	\\
	\label{eq::cr::b}
	\frac{\partial u_\ia}{\partial y}+\frac{\partial v_\ia}{\partial x}&=&
		\sum_{\ib=1}^\dimCMan\bigl[
			(b_{\ia\ib}-e_{\ia\ib})u_\ib-
			(a_{\ia\ib}+f_{\ia\ib})v_\ib
		\bigr].
\end{eqnarray}

If we set $w_\ia=u_\ia+\sqrt{-1}v_\ia$ we easily obtain that the
functions $w_1,\ldots,w_\dimCMan$ satisfy
\begin{eqnarray}
	\label{eq::crcomplex}
	\frac{\partial w_\ia}{\partial\bar z}&=&
		\sum_{\ib=1}^\dimCMan\bigl[
			A_{\ia\ib}w_\ib+
			B_{\ia\ib}\bar{w}_\ib
		\bigr],
		\quad\ia=1,\ldots,\dimCMan,
\end{eqnarray}
where $A_{\ia\ib}$ and $B_{\ia\ib}$ are complex functions of class $C^1$.
By Theorem \ref{thm::ZeroesHyperanalitic} it follows that
the common zeroes of the functions $w_1,\ldots,w_\dimCMan$,
that is $\Leaf\cap\CLocus\cap U$ either coincides with $\Leaf\cap U$
or is a discrete set in $\Leaf\cap U$.

\qed

The following theorem is an easy consequence of Theorem \ref{thm:main}.

\begin{theorem}\label{thm::maintwo}
Let $(\VField,\Eq)$ be a calibrated foliation of class $C^2$
on the complex manifold $\CMan$ with contact locus $\CLocus$.

Let $\ConstA\in\R$ and let $\HyperSurface\sset\CMan$ be the set of
points $\Point\in\CMan$ which satisfies $\Eq(\Point)=\ConstA$.

Let $D\sset\CMan$ be the open set which is the union of all
the leaf $\Leaf$ of $(\VField,\Eq)$ such that
$\Leaf\cap\HyperSurface\neq\void$.

If $\HyperSurface\sset\CLocus$ then $D\sset\CLocus$.
\end{theorem}

\proof
Indeed if $\Leaf$ is a leaf of $(\VField,\Eq)$ such that
$\Leaf\cap\HyperSurface\neq\void$ then 
by hypotesis $\Leaf\cap\HyperSurface\sset\CLocus$.
Since $\Leaf\cap\HyperSurface$ is an integral curve
of the vector field $\VField$ it follows that $\Leaf\cap\HyperSurface$
is not a discrete subset of $\Leaf$ and hence Theorem \ref{thm:main}
implies that $\Leaf\sset\CLocus$.

\qed

\section{\label{section:SectionCauchy}A Cauchy problem}

Let $(\VField,\Eq)$ be a calibrated foliation of class $C^\Smooth$
on the complex manifold $\CMan$.
Then it follows from the definitions that if $\HyperSurface=\{\Eq=\ConstA\}$
is not empty then $\HyperSurface$ is a real hypersurface of class $C^\Smooth$,
and for each $\Point\in\HyperSurface$ we have
$\VField(\Point)\in T_\Point(\CMan)$ and
$J\VField(\Point)\notin T_\Point(\CMan)$.
Moreover we have:

\begin{proposition}\label{prop::AnalyticIntegralCurves}
Let ($\VField$,\Eq) be a calibrated foliation of class $C^\Smooth$, $\Smooth\geq1$
on the complex manifold $\CMan$.
Then the integral curves of $\VField$ are real analytic maps.
\end{proposition}

\proof
Let $g_t:\CMan\to\CMan$ the one parameter group of local transformations
associated to the vector field $\VField$.
Let $p\in\CMan$.
Since $g_{t+s}(\Point)=g_t\bigl(g_s(\Point)\bigr)$ it suffices to prove that
for each $\Point\in\CMan$ the map $t\mapsto g_t(\Point)$ is real analytic
in a neoghbourhood of $t=0$.

Let $\Leaf$ be the leaf passing througth $\Point$ and let $z$
be an adapted holorphic coordinate in a neighbourhood of $p$ in $\Leaf$.
Then, by construction, $z\bigl(g_t(\Point)\bigr)=t$, and hence the map
$t\mapsto g_t(\Point)\in\Leaf$ is real analytic.
Since the inclusion $\Leaf\sset\CMan$ is holomorphic the assertion
easily follows.
\qed

Conversely we have:

\begin{theorem}\label{thm::CauchyCF}
Let $\CMan$ be a complex manifold of complex dimension $\dimCMan$.
Let $\HyperSurface\sset\CMan$ be a closed real hypersurfaces
of class $C^\Smooth$, $\Smooth\geq2$ with holomorphic tangent bundle 
$HT(\HyperSurface)\sset T(\HyperSurface)$.

Let $\VField_0$ be a $C^\Smooth$ vector field on $\HyperSurface$
and let $g^0_t:\HyperSurface\to\HyperSurface$ be the one parameter
group of local transformations associated to the vector field $\VField_0$.
Assume that for each $\Point\in\HyperSurface$ we have
$J\VField(p)\notin T_\Point(\CMan)$ and the map $t\mapsto j\bigl(g^0_t(p)\bigr)$
is real analytic.

Then there exists a neighbourhood $\NeighbA$ of $\HyperSurface$ in $\CMan$
and a calibrated  foliation $(\VField,\Eq)$ of class $C^{\Smooth}$ on $\NeighbA$
such that $\Eq_{|\HyperSurface}=0$, $\VField_{|\HyperSurface}=\VField_0$
and for each leaf $\Leaf$ of $(\VField,\Eq)$ we have $\Leaf\cap\HyperSurface\neq\void$

Such a calibrated foliation is locally unique, that is if 
$(\VField_1,\Eq_1)$ is calibrated  foliation on a neighbourhood $\NeighbB$
of $\HyperSurface$ in $\CMan$ such that
$\Eq_1{|\HyperSurface}=0$ and $\VField_1{|\HyperSurface}=\VField_0$
then there exists a neighbourhood $\NeighbC$ of $\HyperSurface$ in $\CMan$
contained in $\NeighbA\cap\NeighbB$ such that
$\VField=\VField_1$ and $\Eq=\Eq_1$ in $\NeighbC$.

Moreover $(\VField,\Eq)$ is a contact calibrated  foliation on $\NeighbA$,
that is $L_\VField \dc\Eq$ vanishes identically on $\NeighbA$
if, and only if, $\VField_0$ is an infinitesimal symmetry of the distribution
$HT(\HyperSurface)$, that is $[\VField_0,HT(\HyperSurface)]\sset HT(\HyperSurface)$.
\end{theorem}

\proof {
\def\CrossP{W}
\def\CrossMap{\varphi}
\def\CProj{\pi}
Let denote by $\Inclusion:\HyperSurface\to\CMan$ the inclusion map.

By the hypotheses there exists a map
$\CrossP\ni(\Point,z)\mapsto=g_z(\Point)\in\CMan$,
where $\CrossP\sset\HyperSurface\times\C$
is an open subset of $\HyperSurface\times\C$
containing $\HyperSurface\times\{0\}$
such that for each $\Point\in\HyperSurface$ the set
$\CrossP_\Point=\bigl\{z\in\C\mid(\Point,z)\in\CrossP\bigr\}$
is an open connected neighbourhood of $0$ in $\C$,
the function $\CrossP_\Point\ni z\mapsto g_z(p)\in\CMan$ is holomorphic
and if $(\Point,z)\in\CrossP$ with $z=t\in\R$ then
$g_t(\Point)=\Inclusion\bigl(g^0_t(\Point)\bigr)$.
Of course when $t,s\in\R$,
denotig by $\Inclusion_*:T(\HyperSurface)\to T(\CMan)$ the differential
of the inclusion map $\Inclusion:\HyperSurface\to\CMan$, we have
\begin{eqnarray}
\label{eq::generatora::nonumber}
	\left.\frac{d}{dt}g_t(\Point)\right|_{t=0}
	&=&\Inclusion_*\bigl(\VField_0(\Point)\bigr)
\nonumber\\
\label{eq::generatorb::nonumber}
	\left.\frac{d}{ds}g_{is}(\Point)\right|_{s=0}
	&=&J\Inclusion_*\bigl(\VField_0(\Point)\bigr).
\nonumber
\end{eqnarray}

Now set $\CrossP_0=\CrossP\cap\HyperSurface\times\R$
and define $\CrossMap:\CrossP_0\CMan$ putting
$\CrossMap(p,s)=g_{is}(p)$.

Being $\CrossMap(p,0))=\Inclusion(p)$ for each $p\in\HyperSurface$ and
$\dd\CrossMap(p,s)\left(\frac{\partial}{\partial s}\right)=J\VField\bigl(\Inclusion(p)\bigr)$
it follow that after shrinking $\CrossP$ if necessary the map
$\CrossMap:\CrossP_0\to\CMan$ is a diffeomorphism
between $\CrossP_0$ and an open subset $\CrossMap(\CrossP_0)=\NeighbA$ of $\CMan$.

Denoting $\CProj:\CrossP_0\to\R$, $\CProj(\Point,s)=s$ the canonical projection
we set $\Eq=-\CProj\circ\CrossMap^{-1}:\NeighbA\to\R$.
Then $\Eq$ is by construction a function of class $C^\Smooth$ with non vanishing
differential anywhere on $\NeighbA$.

Moreover the formula
\begin{eqnarray}\label{eq::xxx::aa}
	&&G_t\bigl(g_{is}(p)\bigr)=g_{is}\bigl(g_t(p)\bigr)
\nonumber
\end{eqnarray}
defines an one parameter group of local diffeomorphisms of $\NeighbA$.

Let $\VField$ be the infinitesimal generator of $G_t$.
We shall prove that $(\VField,\Eq)$ is a calibrated foliation with
the required properties.

Observe that $\Eq$ is characterized by
\begin{eqnarray}\label{eq::xxx::a}
	&&\Eq\bigl(g_{is}(\Point)\bigr)=-s
\nonumber
\end{eqnarray}
for each $\Point\in\HyperSurface$ and for each $s$ small enought.
Setting $s=0$ we see that $\Eq_{|\HyperSurface}=0$.

We also have
\begin{eqnarray}\label{eq::xxx::x}
	&&G_t\bigl(\Inclusion(p)\bigr)=G_t\bigl(g_0(p)\bigr)=g_0\bigl(g_t(p)\bigr)
	=g_t(p)=\Inclusion\bigl(g^0_t(p)\bigr),
\nonumber
\end{eqnarray}
and hence, for each $\Point\in\HyperSurface$,
\begin{eqnarray}\label{eq::xxx::y}
	&&\VField\bigl(\Inclusion(\Point)\bigr)
		=\left.\frac{d}{dt}G_t\bigl(\Inclusion(p)\bigr)\right|_{t=0}=
		\left.\frac{d}{dt}\Inclusion\bigl(g^0_t(p)\bigl)\right|_{t=0}
		=j_\ast(p)\bigl(\VField_0(\Point)\bigr).
\nonumber
\end{eqnarray}

From
\begin{eqnarray}\label{eq::xxx::nonumbera}
	&&\Eq\Bigl(G_t\bigl(g_{is}(\Point)\bigr)\Bigr)
	=\Eq\Bigl(g_{is}\bigl(g_t(\Point)\bigl)\Bigr)
	=-s=\Eq\bigl(g_{is}(\Point)\bigr)
\nonumber
\end{eqnarray}
we see that the hypersurfaces $\{\Eq=\ConstA\}$ are $G_t$-invariant
and hence $\VField(\Eq)=\dd\Eq(\VField)=0$.

We now prove that $\dc\Eq(\VField)=1$.

We first show that given $\Point\in\VField$, for $z\in\C$, $t\in\R$, 
with $\abs{t},\abs{z}$ small enought
\begin{eqnarray}\label{eq::xxx::nonumberc}
	&&g_{z}\bigl(g_t(\Point)\bigr)=g_{z+t}(\Point).
\nonumber
\end{eqnarray}
Indeed, both sides of (\ref{eq::xxx::nonumberc}) for $p$ and $t$ fixed
are holomorphic functions of $z$.
Since they coincide when $z\in\R$ then they coincide by analytic continuation.

It follows then that
\begin{eqnarray}\label{eq::xxx::nonumberd}
	&&\Eq\bigl(g_{z}(\Point)\bigr)=-{\sf Im} z.
\nonumber
\end{eqnarray}
Indeed, if $z=t+is$ then $s={\sf Im} z$ and
\begin{eqnarray}\label{eq::xxx::nonumberh}
	&&\Eq\bigl(g_{z}(\Point)\bigr)=
	\Eq\Bigl(g_{is}\bigl(g_t(\Point)\bigr)\Bigl)=-s=-{\sf Im} z.
\nonumber
\end{eqnarray}

The formula
$H_t\bigl(g_{is}(\Point)\bigr)=g_{i(s-t)}(\Point)$
defines an one parameter group of local diffeomorphisms of $\NeighbA$.
We now show that the infinitesimal generator of $H_t$ is $-J\VField$.
Indeed we have
\begin{eqnarray}\label{eq::xxx::nonumberb}
	\left.\frac{d}{dt}H_t\bigl(g_{is}(\Point)\bigr)\right|_{t=0}
	&=&\left.\frac{d}{dt}g_{i(s-t)}(\Point)\right|_{t=0}
	=-J\left.\frac{d}{dt}g_{(t+is)}(\Point)\right|_{t=0}
\nonumber\\
	&=&-J\left.\frac{d}{dt}g_{is}\bigl(g_t(\Point)\bigr)\right|_{t=0}
	=-J\left.\frac{d}{dt}G_t\bigl(g_{is}(\Point)\bigr)\right|_{t=0}
\nonumber\\
	&=&-J\VField(\Point).
\nonumber
\end{eqnarray}

Thus we obtain
\begin{eqnarray}\label{eq::xxx::nonumbere}
	\Eq\Bigl(H_t\bigl(g_{is}(\Point)\bigr)\Bigr)
	=\Eq\bigl(g_{i(s-t)}(\Point)\bigr)=t-s,
\nonumber
\end{eqnarray}
and hence
\begin{eqnarray}
\label{eq::xxx::nonumberg}
	\dc\Eq(\VField)\bigl(g_{is}(\Point)\bigr)
	&=&-\dd\Eq(J\VField)\bigl(g_{is}(\Point)\bigr)
\nonumber\\
\label{eq::xxx::nonumberf}
	&=&\left.\frac{d}{dt}\Eq\Bigl(H_t\bigl(g_{is}(\Point)\bigr)\Bigr)\right|_{t=0}
	=\left.\frac{d(t-s)}{dt}\right|_{t=0}
	=1.
\nonumber
\end{eqnarray}

We end the proof that $(\VField,\Eq)$ is a calibrated foliation
showing that $[\VField,J\VField]=0$.

It suffices to prove that $G_t$ and $H_t$ commute.
Indeed we have
\begin{eqnarray}
	H_{t_1}\circ G_{t_2}\bigl(g_{is}(p)\bigr)=
	G_{t_2}\circ H_{t_1}\bigl(g_{is}(p)\bigr)=
	g_{i(s-t_1)}\bigl(g_{t_2}(p)\bigr),
\nonumber
\end{eqnarray}
and the proof of the existence of a calibrated foliation is completed.

We now prove the uniqueness of $(\VField,\Eq)$.
Let $(\VField_1,\Eq_1)$ an other calibrated  foliation on
a neighbourhood $\NeighbB$ of $\HyperSurface$ in $\CMan$ such that
$\Eq_1{|\HyperSurface}=0$ and $\VField_1{|\HyperSurface}=\VField_0$.

The leafs of the foliation $(\VField,\Eq)$ and the ones of $(\VField_1,\Eq_1)$
which intersect $\HyperSurface$ both intersect $\HyperSurface$ along the
integral curves of $\VField_0$ and therefore are the same.

Let hence $\NeighbC$ be the union of all the leafs of the foliation $(\VField,\Eq)$
which intersect $\HyperSurface$.
Then the restriction of the function $w=\Eq-\Eq_1$ to each leaf $\Leaf$
is an harmonic function on $\Leaf$ which vanishes $\Leaf\cap\HyperSurface$.
Since $\dc w(\Point)(\VField_0)=
\dc\Eq(\Point)(\VField_0)-\dc\Eq_1(\Point)(\VField_0)=1-1=0$
for each $\Point\in\HyperSurface$
it follows that $\dc w_{|\Leaf\cap\HyperSurface}=0$.
By lemma \ref{lemma::HarmonicZero} below it follows that $w_{|\Leaf}=0$,
that is $\Eq_{|\Leaf}={\Eq_1}_{|\Leaf}$.
Since the leaf $\Leaf\sset\NeighbC$ is arbitrary then
$\Eq_{|\NeighbC}={\Eq_1}_{|\NeighbC}$
and also $\VField_{|\NeighbC}={\VField_1}_{|\NeighbC}$
easily follows.

It remains to prove the last assertion of the Theorem.

Let $\omega=L_\VField \dc\Eq$.
We have to prove that $\omega$ vanishes identically on $\NeighbA$
if, and only if, $[\VField_0,HT(\HyperSurface)]\sset HT(\HyperSurface)$.

By Theorem \ref{thm::maintwo} it suffices to prove that
$[\VField_0,HT(\HyperSurface)]\sset HT(\HyperSurface)$
if and only if $\omega$ vanishes identically on $\HyperSurface$.

Let $\Point_0\in\HyperSurface$.
Let $U$ be a neighbourhood of $\Point_0$ in $\CMan$ and vector fields
$X_1,\ldots,X_\dimCMan$ such that
$\VField,J\VField,X_1,JX_1\ldots,X_\dimCMan,JX_\dimCMan$
is a frame for $T(\CMan)$ on $U$ such that
$\VField,X_1,JX_1\ldots,X_\dimCMan,JX_\dimCMan$
is a frame for ${\sf Ker}\, \dd\Eq$ and
$X_1,JX_1\ldots,X_\dimCMan,JX_\dimCMan$
is a frame for ${\sf Ker}\, \dd\Eq\cap{\sf Ker}\,  \dc\Eq$.

As in the proof of Theorem \ref{thm:main} we see that
that $\omega$ vanishes identically on $\HyperSurface$
if, and only if for $\ia=1,\ldots,\dimCMan$ the functions
$\dc\Eq\bigl([\VField,X_\ia]\bigr)$ and $\dc\Eq\bigl([\VField,JX_\ia]\bigr)$
vanish identically on $\HyperSurface$, that is,
since $T(\HyperSurface)={\sf Ker}\, \dd\Eq_{|\HyperSurface}$ and
$HT(\HyperSurface)={\sf Ker}\, \dd\Eq_{|\HyperSurface}\cap{\sf Ker}\,  \dc\Eq_{|\HyperSurface}$,
if, and only if for $\ia=1,\ldots,\dimCMan$
we have
$[\VField,X_\ia]_{|\HyperSurface}\in HT(\HyperSurface)$ and
$[\VField,JX_\ia]_{|\HyperSurface}\in HT(\HyperSurface)$.

But we have
$$
	[\VField,X_\ia]_{|\HyperSurface}=[\VField_{|\HyperSurface},{X_\ia}_{|\HyperSurface}]
	=[\VField_0,{X_\ia}_{|\HyperSurface}]
$$
and
$$
	[\VField,JX_\ia]_{|\HyperSurface}=[\VField_{|\HyperSurface},{JX_\ia}_{|\HyperSurface}]
	=[\VField_0,{JX_\ia}_{|\HyperSurface}]
$$

Being ${X_1}_{|\HyperSurface},{JX_1}_{|\HyperSurface}\ldots,
{X_\dimCMan}_{|\HyperSurface},{JX_\dimCMan}_{|\HyperSurface}$
be a frame for $HT(\HyperSurface)$ it follows that
$\omega$ vanishes identically on $\HyperSurface$
if, and only if, $[\VField_0,HT(\HyperSurface)]\sset HT(\HyperSurface)$.

The proof of the Theorem is therefore completed.

\qed }

The idea to construct the function $\Eq$ as the imaginary part of the
function obtained by complex analytic continuation of the integral
curves of a vector field is taken from \cite{article:DuchampKalka}.

\begin{lemma}\label{lemma::HarmonicZero}
Let $D\sset\C$ be a domain such that $D\cap\R\neq\void$ and let
$w:D\to\R$ be a harmonic function. If $w$, $\partial w/\partial x$ and $\partial w/\partial y$ vanish on $D\cap\R$ then $w$ vanishes identically on $D$.
\end{lemma}

\proof
Let $x_0\in D\cap\R$. Then there exists a convex neighbourhood $U\sset D$ of $x_0$ and
an holomorphic function $f:U\to\C$ such that $w={\sf Re} f$.
By the Cauchy-Riemann equations it follows that $f'(x)=0$ on the interval $U\cap\R$.
By the analytic continuation principle for the holomorphic functions it follows that
$f'(z)=0$ on $U$ and hence $f(z)$ is constant on $U$
But then  $w={\sf Re} f$ also is constant on $U$. Since $w$ vanishes on $U\cap\R$ then
it vanishes also on $U$.
Since $w$ is real analytic then it must vanish identically on $D$.
\qed

\section{\label{section:ComplexMA}The complex Monge-Amp\`ere equations}

Let $\CMan$ be a connected complex manifold of complex dimension $\dimCMan+1$.

Let $\Eq\in C^\Smooth(\CMan)$, $\Smooth\geq2$ be a function without critical point.
For each constant $c$ let denote by $\HyperSurface_\ConstA$ the set of points of $M$ where
$\Eq$ assume the value $\ConstA$.

We denote by $H_\Eq\sset T(\CMan)$ the distribution ${\sf Ker}\, \dd\Eq\cap{\sf Ker}\,  \dc\Eq$.
Assume that $\dd\Eq\wedge \dc\Eq\wedge(\dd\dc\Eq)^{\dimCMan}$ do not
vanish on $\CMan$. Then each hypersuface $\HyperSurface_c$
is a ${\rm CR} $ contact manifold with $HT(\HyperSurface)={H_\Eq}_{|\HyperSurface}$
and contact form $\dc\Eq_{|T(\HyperSurface)}$.
We then denote by $\VField_\Eq$ the unique vector field on $\CMan$
(of class at least $C^{\Smooth-2}$) 
which satisfies $\dd\Eq(\VField_\Eq)=0$, $\dc\Eq(\VField_\Eq)=1$ and
$\dd\dc\Eq(\VField_\Eq,X)=0$ for each vector field $X$ which satisfies $\dd\Eq(X)=0$.

In other word $\VField_\Eq$ is the vector field on $\CMan$ which is tangent to
each hypersurface $\HyperSurface_\ConstA$ and coincides on $\HyperSurface_\ConstA$
with the Reeb vector field associated to the contact form
$\dc\Eq_{|T(\HyperSurface)_\ConstA}$.
Observe that $\VField_\Eq$ can be characterizez by the conditions
$\dd\Eq(\VField_\Eq)=0$, $\dc\Eq(\VField_\Eq)=1$ and
$[\VField_\Eq,H_\Eq]\sset H_\Eq$.


\begin{lemma}\label{lemma::contactversusma}
Let $(\VField,\Eq)$ be a calibrated foliation of class $C^2$
on the complex manifold $\CMan$ of complex dimension $\dimCMan+1$
with contact locus $\CLocus$.

Then we have the identity
\begin{eqnarray}\label{eq::LieMA}
	&&\VField\res\dd\dc\Eq=L_\VField\dc\Eq
\end{eqnarray}
and the form $(\dd\dc\Eq)^{\dimCMan+1}$ vanishes on $\CLocus$.
\end{lemma}

\proof
Let $X$ be a vector field on $\CMan$.
Since $X\bigl(\dc\Eq(\VField)\bigr)=X(1)=0$ then
\begin{eqnarray}
\label{eq::LieRelation::nonumber}
	\dd\dc\Eq(\VField, X)
	&=&\VField\bigl(\dc\Eq(X)\bigr)-X\bigl(\dc\Eq(\VField)\bigr)
		-\dc\Eq\bigl([\VField,X]\bigr)
\nonumber\\
	&=&L_\VField \dc\Eq(X)-X\bigl(\dc\Eq(\VField)\bigr)
	=L_\VField \dc\Eq(X).
\nonumber
\end{eqnarray}

It follows that if $\Point\in\CLocus$
then $\VField(\Point)\neq0$ belongs to the radical of the bilinear form
$(X,Y)\mapsto \dd\dc\Eq(X,JY)$, and hence its rank is strictly less than
$\dimCMan+1$ and this implies that $(\dd\dc\Eq)^{\dimCMan+1}$ vanishes at
$\Point$.

\qed

\begin{theorem}\label{thm::CauchyMA}
Let $\CMan$ be a connected complex manifold of complex dimension $\dimCMan+1$
Let $\HyperSurface\sset\CMan$ be a closed contact CR hypersurface of class $C^\Smooth$,
$\Smooth\geq3$ and let $\CForm$ be a contact form for $(\HyperSurface, HT(\HyperSurface))$
of class $C^{\Smooth-1}$ with associated Reeb vector field $\VField_0$.
Let denote by $\Inclusion:\HyperSurface\to\CMan$ the inclusion map.

Assume that the Reeb vector field $\VField_0$ is of class $C^\Smooth$.

Then the following conditions are equivalent:

\begin{enumerate}
	\item\label{item::AnalyticFlow}
	for each integral curve $\gamma(t)$ of $\VField_0$ the map
	$t\mapsto\Inclusion\bigl(\gamma(t)\bigr)$ is real analytic;
	\item\label{item::MASolution}
	there exists a open neighbourhood $\NeighbA\sset\CMan$ of
	$\HyperSurface$ and a function $\Eq\in C^\Smooth(\NeighbA)$ which satisfies
	\begin{equation}\label{eq::CauchyMongeAmpere}
		\left\{
		\begin{array}{l}
		(\dd\dc\Eq)^{\dimCMan+1}=0\ {\rm on}\ \NeighbA,\\
		\\
		\dd\Eq\wedge \dc\Eq\wedge(\dd\dc\Eq)^{\dimCMan}\neq0\ {\rm on}\ \NeighbA,\\
		\\
		\Eq_{|\HyperSurface}=0,\\
		\\
		\dc\Eq_{|T(\HyperSurface)}=\CForm.
		\end{array}
		\right.
	\end{equation}
\end{enumerate}
\end{theorem}

\proof
Assume (\ref{item::AnalyticFlow}).
By theorem \ref{thm::CauchyCF} there exists a neighbourhood $\NeighbA\sset\CMan$
of $\HyperSurface$ and a contact calibrated foliation $(\VField,\Eq)$ on $\NeighbA$
such that $\VField_{|\HyperSurface}=\VField_0$.
We claim that the function $\Eq$ satisfies (\ref{eq::CauchyMongeAmpere})
By lemma \ref{lemma::contactversusma} the function $\Eq$ satisfies the
Monge-Amp\`ere equation.

By construction $\Eq_{|\HyperSurface}=0$ and for each $\Point\in\HyperSurface$
we have ${\sf Ker}\,  \dc\Eq_{|T_\Point(\HyperSurface)}={\sf Ker}\, \CForm(\Point)$.
Since
$$
\dc\Eq_{|T(\HyperSurface)}(\VField_0)=\dc\Eq(\VField)_{\HyperSurface}=1=\CForm(\VField_0),
$$
it follows that $\dc\Eq_{|T(\HyperSurface)}=\CForm$.

Finally, since $\HyperSurface$ and a contact calibrated foliation then 
$\CForm\wedge(\dd\CForm)^\dimCMan$ does not vanish on $\HyperSurface$ and
it is then easy to show that shrinking the neighbpurhood $\NeighbA$ if necessary,
the function $\Eq$ satisfies $\dd\Eq\wedge \dc\Eq\wedge(\dd\dc\Eq)^{\dimCMan}\neq0$
in a neighbourhood of $\HyperSurface$ in $\CMan$.

Conversely assume that $\Eq$ is a solution of class $C^\Smooth$, $\Smooth\geq3$
of (\ref{eq::CauchyMongeAmpere}) in a neighbourhood $\NeighbA$
of $\HyperSurface$ in $\CMan$.

Set $\VField=\VField_\Eq$. Then $\VField$ is a $C^{\Smooth-2}$ vector field
on $\CMan$ which at each point $\Point\in\HyperSurface$ is tangent to
$\HyperSurface$ and on coincides with the Reeb vector field $\VField_0$.

It suffices then to prove that the integral curves of the vector field $\VField$
are analytic curves in $\CMan$.

Let $\Point$ be an arbitrary point of $\CMan$
Let $U$ be a neighbourhood of $\Point$ in $\CMan$ and
let $X_1,\ldots,X_\dimCMan$ be $C^\Smooth$ vector fields such that
$X_1,JX_1\ldots,X_\dimCMan,JX_\dimCMan$ is a local frame on $U$ for
the distribution $H_|Eq$.

Then $\VField,J\VField, X_1,JX_1\ldots,X_\dimCMan,JX_\dimCMan$
is a local frame on $U$ for the tangent bundle  $T(\CMan)$.

By construction $\VField$ is a $C^{\Smooth-2}$ vector field which satisfies
$\dd\Eq(\VField)=0$, $\dc\Eq(\VField)=1$ and
$\dd\dc\Eq(\VField,X_\ia)=\dd\dc\Eq(\VField,JX_\ia)=0$
for $\ia=1,\ldots,\dimCMan$.

Since $(\dd\dc\Eq)^{\dimCMan+1}=0$ then
\bigskip
$$
0=\det \left(\begin{array}{ccccc}
\dd\dc\Eq(\VField,J\VField) &\dd\dc\Eq(\VField,X_1)
 & \dots   &\dd\dc\Eq(\VField,X_n)
\\
\dd\dc\Eq(\VField,X_1)
 & \dd\dc\Eq(X_1,JX_1) & \dots   & \dd\dc\Eq(X_1,JX_\dimCMan) \\
\vdots & \vdots &  &   \vdots \\
\dd\dc\Eq(\VField,JX_\dimCMan)&\dd\dc\Eq(X_1,JX_\dimCMan)& \dots   &\dd\dc\Eq(X_\dimCMan,JX_\dimCMan)  
\end{array}\right)
$$
\bigskip
$$
=\det \left(\begin{array}{ccccc}
\dd\dc\Eq(\VField,J\VField) &0
 & \dots   &0
\\
0& \dd\dc\Eq(X_1,JX_1) & \dots   & \dd\dc\Eq(X_1,JX_\dimCMan) \\
\vdots & \vdots &  &   \vdots \\
0&\dd\dc\Eq(X_1,JX_\dimCMan)& \dots   &\dd\dc\Eq(X_\dimCMan,JX_\dimCMan)  
\end{array}\right)
$$
\bigskip
$$
=\dd\dc\Eq(\VField,J\VField)
\det \left(\begin{array}{ccccc}
\dd\dc\Eq(X_1,JX_1) & \dots   & \dd\dc\Eq(X_1,JX_\dimCMan) \\
\vdots & \vdots &\\
\dd\dc\Eq(X_1,JX_\dimCMan)& \dots   &\dd\dc\Eq(X_\dimCMan,JX_\dimCMan)  
\end{array}\right).
$$
The last determinant does not vanish so 
$$
\dd\dc\Eq(\VField,J\VField)=0
$$
and since clearly $\dd\dc\Eq(\VField,\VField)=0$ we obtain $\VField\res\dd\dc\Eq=0$ on $U$.

Since the point $\Point\in\CMan$ is arbitrary it follows that
$\VField\res\dd\dc\Eq=0$ on $\CMan$.

For each vector field $X$ we also have
$$
	\dd\dc\Eq(J\VField,X)=-\dd\dc\Eq(\VField,JX)=0
$$
that is $J\VField\res\dd\dc\Eq=0$.

Let us prove that $[\VField,J\VField]=0$.
By Theorem 2.4, pag. 549 of \cite{article:BedfordKalka}
it follows in particular that
for each $\Point\in\CMan$ the the vector subspace of $T_\Point(\CMan)$
of the vectors $Z\in T_\Point(\CMan)$ which satisfy
$Z\res\dd\dc\Eq=0$ is a $J$-invarant subspace of real dimension 2.

Hence we have $[\VField,J\VField]=a\VField+bJ\VField$ for some functions $a$ and $b$.

We prove that the functions $a$ and $b$ vanish on $\CMan$.
By (\ref{eq::lemma::c}) of Lemma \ref{lemma::ElementaryLemma}
it follows that $\dd\Eq\bigl([\VField,J\VField]\bigr)=0$ and hence
$$0=\dd\Eq\bigl([\VField,J\VField]\bigr)=a\dd\Eq(\VField)+b\dd\Eq(J\VField)=-b.$$

By \ref{eq::lemma::b} of Lemma \ref{lemma::ElementaryLemma}
it follows that $\dc\Eq\bigl([\VField,J\VField]\bigr)=-\dd\dc\Eq(\VField,J\VField)$
and hence
$$0=\dc\Eq\bigl([\VField,J\VField]\bigr)=a\dc\Eq(\VField)+b\dc\dc\Eq(J\VField)=a.$$

Thus we have proved that $(\VField,\Eq)$ is a calibrated foliation of class at least $C^1$.
The analiticity of the integral curves of $\VField$ follows then from
Proposition \ref{prop::AnalyticIntegralCurves}.
\qed

\begin{remark}
When $\HyperSurface$ and $\CForm$ are real analytic the
previous theorem is an immediate consequence of
Proposition 1.5 of  \cite{article:BedfordBurns}.
Proposition 1.1 of \cite{article:BedfordBurns} gives
the uniqueness of a $C^3$ solution of the problem (\ref{eq::CauchyMongeAmpere}).
\end{remark}

\end{document}